\documentclass[12pt,a4paper,notitlepage]{article}
\usepackage{t1enc}
\usepackage[latin2]{inputenc}

\usepackage{amsthm}
\usepackage[intlimits]{amsmath}
\usepackage{amssymb}
\usepackage{amsthm}
\usepackage{graphics}

\title{Random tree growth with general weight function}
\author{Anna Rudas
\\
{Institute of  Mathematics}
\\
{Technical University Budapest}
}

\date{October 22, 2004}

\newtheorem*{thm}{Theorem}
\newtheorem*{lemma}{Lemma}
\theoremstyle{definition} \newtheorem*{prf}{Proof}

\newcommand{\E}{\mathbb{E}}

\newcommand{\N}{\mathbb{N}}
\newcommand{\R}{\mathbb{R}}
\newcommand{\pr}{\mathbb{P}}
\newcommand{\de}{\text{deg}}
\newcommand{\ind}{1\!\!1}
\newcommand{\lev}{\bar}

\begin{document}

\maketitle

\begin{abstract}
We extend the results of B.~Bollob\'as, O.~Riordan, J.~Spencer, G.~Tusn\'ady 
\cite{Bollobas_degree} and M\'ori \cite{Mori}. 
We consider a  model of random tree growth, where
at each time unit a new node is added and attached to an already
existing node chosen at random. The probability with which a node
with degree $k$ is chosen is  proportional to $w(k)$, where
$w:\N\to\R_+$ is a fixed weight function. 
We prove that if $w$ fulfills some asymptotic requirements
then  the decay of the  distribution of the degree sequence as
$k\to\infty$    is $\text{exp}\{-c\sum_{j=0}^{k-1}w(j)^{-1}\}$,
where $c$ is an  explicitly  computable positive constant. In
particular, if $w(k)$ is asymptotically linear then this decay
has  power law. 
Our method of proof is analytic rather than combinatorial, having
the advantage of robustness: only asymptotic properties of  the
weight function $w(k)$ are used, while in \cite{Bollobas_degree}
and \cite{Mori} the explicit law $w(k)=ak+b$ is assumed. 
   

\end{abstract}

\section{Introduction}

In 1999, Barab\'asi and Albert \cite{Barabasi} suggested a random graph process for modeling the growth 
of real-world networks. In this model each time a new node is added to the system, it is attached to a fixed 
number  of already existing nodes, that are selected with probabilities proportional to their degrees. Barab\'asi 
and Albert argued that in this growth model, after many steps, the degree sequence obeys a power law decay
(the proportion $P(k)$ of nodes with degree exactly $k$ decays as $k^{-\gamma}$). According to 
measurements, this scale-free property seems to be true for many real-world networks, such as the world wide 
web. 

In 2000, Bollob\'as, 
Riordan, Spencer and Tusn\'ady \cite{Bollobas_degree} gave a mathematically exact definition of the model 
described above, and obtained $P(k)$ asimptotically uniformly for all $k\le n^{1/15}$, where $n>>1$ is the 
number of 
nodes, proving as a consequence that $\gamma=3$. Their proof depends on a combinatorial method that 
relates the tree model to random pairings, and also on the Azuma-Hoeffding inequality. Independently, in 
\cite{Mori}, M\'ori proved similar results for trees,
also with martingale techniques, although different from those used in \cite{Bollobas_degree}. 
For a survey on mathematical results on scale-free random graphs, see \cite{Bollobas_summary}.

Our model generalizes the rule of preferential attachment: we let a weight function $w: \N \to \R ^+$ 
determine the law with which a node is selected. A node with degree $k$ is chosen with probability 
proportional to $w(k)$. In the proofs of both \cite{Bollobas_degree} and \cite{Mori}, it is essential that 
the weight function has the exact form $w(k)=ak+b$. We treat general weight functions, with some restrictions 
on the asymptotics of $w$.
   
We define the model for trees, in each step the new node attaches to one already existing node. 
We prove that for any $k$, the proportion of nodes with at 
least $k$ outgoing degree converges in probability to $c_k$, a constant that can be computed, given $w$. The 
decay of $c_k$ as a function of $k$ is given by 
\[
\text{exp}\{-c\sum_{j=0}^{k-1}\frac{1}{w(j)}\},
\]
where $c$ is a computable positive constant.
This means power order decay for any $w$ that is asimptotically linear.

Our method of proof is analytic rather than combinatorial, having the advantage of robustness. We rely only 
on asymptotic properties, and not on the exact form of the weight function $w$.

In \cite{Bollobas_degree} any fixed number $m\ge 1$ can describe the number of edges that a newly appearing 
node sends to the already existing ones, so the result obtained there is not only valid for trees. 
The method they use is to construct the $(m>1)$-model directly from the $(m=1)$-model, 
which procedure fails in the 
case of general weight functions. This is the reason why our result for trees can not
be simply modified to be valid for $m>1$.

\section{Model}

Let a strictly positive weight function $w: \N \to \R^+$ be given with the following restrictions.
We demand that $w(k)\to\infty$ as $k\to \infty$, and that 
the increments of $w$ be bounded from above: there
exists an $r<\infty$ for which
\begin{equation}\label{eq:felt1}
w(k+1)-w(k)\le r
\end{equation} 
for all $k$. (We do not need monotonicity for $w$.) Let us fix $w(0)=1$, it will be apparent 
that this can be done without loss of generality. 
We also demand that the weight function be of the form
\begin{equation}\label{eq:felt2}
w(k)=k^{\alpha}+v(k)
\end{equation}
with some $0< \alpha \le 1$ and $v(k)$ fulfilling the requirements below. 
\newline (\emph{a}) If $\alpha =1$, then let 
\begin{equation}\label{eq:felt_a}
v(k)=o(k) \text{ as }k\to \infty , \text{ and }
\sum _{k=1}^{\infty} \frac{v(k)}{k^2}< \infty,
\end{equation}
(b) If $0< \alpha < 1$, then let
\begin{equation}\label{eq:felt_b}
v(k)=o(k^{\alpha}) \text{ as }k\to \infty.
\end{equation}
We believe that these requirements can be weakened.


\subsection{Discrete time model}

We let a random tree $\text{TREE}_{\text{discr}}(n)$ grow in discrete time 
in the following way. At $n=0$ 
one single node exists with no edges. At each time $n\ge 1$ a new node is attached by an edge 
to one of the already existing nodes, chosen randomly, with a distribution depending on
the degree structure of $\text{TREE}_{\text{discr}}(n-1)$. We denote by 
$\de (j, n)$ the outgoing degree of
node $j$ at time $n$, namely the number of neighbors of node $j$ that have appeared after it.
The probability that node $n$ attaches to node $j$ is 
\[
p(j,n-1)=\frac{w(\de (j, n-1))}{\sum _{i=0}^{n-1}w(\de (i, n-1))}.
\]

From this it is obvious that we can choose $w(0)=1$,
multiplying the weight function with a positive constant does not change the 
process.

\subsection{Continuous time model}

We define another random tree, $\text{TREE}_{\text{cont}}(t)$, growing in continuous time as follows.
Again we speak of the outgoing degree $\de (j,t)$ of node $j$ at time $t$, this means the number of 
those neighbors of node $j$ at time $t$ that have appeared after it. 

At $t=0$ one node, the root, exists with no edges. 
At any time $t\ge 0$ all of the existing nodes, independently of each other, ``give birth to a child'' 
with a rate determined by the actual outgoing degree of the node: a new node appears and is attached to node 
$j$ with rate $\de (j,t)$. This way the probability that the next node that appears in 
$\text{TREE}_{\text{cont}}(t)$ is attached to node $j$ is
\[
p(j,t)=\frac{w(\de (j, t))}{\sum _{i=0}^{N_t}w(\de (i, t))},
\]
where $N_t$ denotes the number of nodes in $\text{TREE}_{\text{cont}}(t)$, apart from the root.

The two models are equivalent in the following sense. For $n\in \N$ let 
$T_n:=\inf \{t: N_t=n\}$ be the random stopping time when node $n$ appears in 
the tree. With these notations
\[\text{TREE}_{\text{cont}}(T_n)\stackrel{\textrm{d}}{=}\text{TREE}_{\text{discr}}(n).\]

\section{Result and sketch of proof}

\begin{thm}
Let 
\[
N_k(t):=\#\{\text{nodes that have at least }k\text{ children in } \text{TREE}_{\text{cont}}(t)\}.
\]
With this notation, for any fix $k\in \N$
\[
\frac{N_k(t)}{N_0(t)}\stackrel{\pr}{\to} c_k,
\]
these constants are given by
\begin{equation}
c_k:=\prod_{j=0}^{k-1}\frac{w(j)}{\lambda ^*+w(j)}. \label{eq:c_def}
\end{equation}
Here $\lambda ^*$ is the unique solution of the following equation: 
\begin{equation}\label{eq:lambdadef}
\sum_{k=1}^{\infty}\prod _{j=0}^{k-1}\frac{w(j)}{\lambda^*+w(j)}=1.
\end{equation}
 
\end{thm}

Note that if the theorem holds, (\ref{eq:lambdadef}) has to be true, since for all $t>0$
$\sum _{j=1}^{\infty}N_k(t)$ simply counts the number of edges apparent in $\text{TREE}_{\text{cont}}(t)$, 
so divided by $N_0(t)$ it must be 1 in the limit.

Note also that in the special case $w(j)=j+1$, which corresponds to the Barab\'asi model,
$\lambda ^*=2$ and so $c_k=\prod_{j=1}^k\frac{j}{2+j}=\frac{2}{(k+1)(k+2)}$.
The ratio of nodes with exactly $k$ outgoing edges is then 
$c_k-c_{k+1}=\frac{4}{(k+1)(k+2)(k+3)}.$ The result for trees in 
\cite{Bollobas_degree} is exactly this.

\subsection{Sketch of proof}

The main advantage of the continuous time model is that when a new node appears, it starts
to give life to a new subgraph that (in distribution) behaves just like the whole tree. 
Certain subgraphs, e.g.~the one below the first child of the root, and the one below the second, grow 
independently of each other. 

Therefore it is essential to study the process $X_t$, that is
the number of children of the root at time $t$. We do this in some preliminary steps.

With the help of these results we 
prove the theorem in the following way. 
We define $\tilde{n}_k(t):=e^{-\lambda ^* t}\E(N_k(t))$ and 
$\tilde{n}_{k,l}(t):=e^{-2\lambda ^* t}\E(N_k(t)N_l(t))$ with $\lambda ^*$ defined by 
(\ref{eq:lambdadef}). The first step is to prove that
\[
\tilde{n}_k(t)\to d_k \text{ as } t\to \infty
\]
where $d_k$ are constants for which $\frac{d_k}{d_0}=c_k$. In the second step we prove that
\[
\tilde{n}_{k,l}(t)\to Cd_kd_l \text{ as } t\to \infty
\]
with some constant $0<C<\infty$. Using these we finally show in step three that the theorem 
holds.

\section{Preliminary steps}\label{chap:propertiesofxt}

\subsection{Properties of $\varrho (s)$}

Let $\tau _j$ be the time of birth of the $j^{\mathrm{th}}$ child of the root, with this 
\[
X_t:=\sum_{j=1}^{\infty} \ind_{ \{ \tau_j<t \} }
\] 
is the number of children of the root at time $t$.
We denote by $\varrho (s)$ the density of the point process $\tau=(\tau_j)_{j=1}^{\infty}$ on 
$\R^+$,
\[
\varrho (s):=\lim_{\varepsilon \to 0}\varepsilon ^{-1}
\pr(\text{a new child of the root appears in }(s,s+\varepsilon)).
\] 

It will be useful later, during the proof of the theorem, to be able to compute the Laplace 
transform of $\varrho (s)$. For any function $\phi$ we denote its Laplace transform 
by $\hat {\phi}$:
\[
\hat{\phi}(\lambda):=\int_0^{\infty} e^{-\lambda t}\phi (t)\mathrm{d}t.
\]

Consider that
\[
\varrho(s)=\lim_{\varepsilon\to 0}\left \{ \sum_{k=1}^{\infty} \varepsilon ^{-1}
\pr(\tau _k \in (s,s+\varepsilon))+\varepsilon ^{-1}o(\varepsilon)\right \}=
\sum_{k=1}^{\infty}\varphi _k(s),
\]
where $\varphi _k(s)$ is the density of the random variable $\tau_k$. Note that $\tau _k$ is 
the sum of independent, 
exponentially distributed random variables, $\tau_k=\sum_{j=1}^{k}\theta _j$ with
distribution $\pr (\theta _j<t)=1-e^{-w(j-1)t}$. This way
we can easily obtain the Laplace transform of $\varphi _k(t)$:
\[
\hat{\varphi}_k(\lambda)=\E(e^{-\lambda \tau_k})=\prod_{j=1}^{k} 
\E(e^{-\lambda\theta _j})=
\prod_{j=0}^{k-1}\frac{w(j)}{\lambda+w(j)}.
\]
This leads to an explicit formula of the Laplace transform of $\varrho$:
\begin{equation}\label{eq:varrhokalap}
\hat{\varrho}(\lambda)=\int_0^{\infty}e^{-\lambda t}\varrho(t)\mathrm{d}t=
\int_0^{\infty}\sum_{k=1}^{\infty}e^{-\lambda t}\varphi _k(t)\mathrm{d}t=
\sum_{k=1}^{\infty}\prod_{j=0}^{k-1}\frac{w(j)}{\lambda+w(j)}.
\end{equation}

It can be easily seen how the finiteness of this sum depends on $\lambda$:  
in the case $0<\alpha <1$, $\hat{\varrho}(\lambda)$ is finite for any $\lambda>0$, 
and in the case $\alpha =1$, it is finite for any $\lambda>1$ 
(and is infinite for $\alpha\le 1$). We will use this fact later.

Recalling the
definition of $\lambda^*$ in (\ref{eq:lambdadef}), note t hat it is equivalent with 
$\hat{\varrho}(\lambda^*)=1$.
The reason for this will be apparent shortly. For later use we introduce the constant
\begin{equation}\label{eq:b1def}
B_1:=\int_{0}^{\infty}e^{-2\lambda^* s}\varrho(s)\mathrm{d}s=\hat{\varrho} (2\lambda^*)<1.
\end{equation}

\subsection{Estimate of $\varrho _2(u,s)$}\label{chap:varrho2}

For $u\ne s$ we define the second correlation function of the point process $\tau$:
\[
\varrho_2(u,s):=\lim_{\varepsilon, \delta \to 0}(\varepsilon \delta)^{-1}
\pr ((u,u+\varepsilon)\text{ and }(s,s+\delta)\text{ both contain a point from $\tau$}).
\]
For $u=s$ we define $\varrho _2(u,u):=0$.

\begin{lemma}
The following integral is finite:
\begin{equation}\label{eq:kell}
B_2:=\int_0^{\infty}\int_0^{\infty}e^{-\lambda ^* (u+s)}\varrho_2(u,s)\mathrm{d}u\mathrm{d}s
<\infty.
\end{equation}
\end{lemma}

\begin{prf}
In order to estimate $\varrho_2(u,s)$, let us examine the Markov process $X_t$. Its 
infinitesimal generator is $G$ with 
\[G_{i,j}=w(i)(-\delta _{i,j}+\delta _{i+1,j}).\]


We can easily determine the eigenfunctions $f_r$ of $G$ with eigenvalue $r$
for arbitrary $r\in (0,\infty )$: we need
\[-w(k)f_r(k)+w(k)f_r(k+1)=rf_r(k)\]
for all $k$. Equivalently
\begin{equation}\label{eq:frdef}
f_r(k)=\prod_{j=0}^{k-1}\frac{r+w(j)}{w(j)}f_r(0)
\end{equation}
(let $f_r(0)=1$.)
With this, $M_t:=e^{-rt}f_r(X_t)$ is a positive martingale with $M_0=1$, so 
\begin{equation}\label{eq:martingale}
\E (f_r(X_t))=e^{rt}. 
\end{equation}

For estimating $\varrho_2(u,s)$, the following will be useful.  
Let us extend $f_r$ piecewise linearly to the positive real numbers, we will refer to the
extension with the same letter. Observe that we can find a finite $r$ for which $f_r$ is convex:
\begin{multline*}
{f}_r(k+1)-2{f}_r(k)+{f}_r(k-1)\\
={f}_r(k)\frac{r}{w(k)(r+w(k-1))}
\left(r-(w(k)-w(k-1))\right),
\end{multline*}
this is positive if $r>w(k)-w(k-1)$. From assumption (\ref{eq:felt1}) we have on $w$, there
exists such an $r<\infty$.

Let us fix such an $r$ and begin the estimate of $\varrho_2(u,s)$. We introduce a generic constant 
$D$. The exact value of $D$ may change from line to line and this will not cause confusion: the only 
thing we will need is that $D$ always denotes a finite number. 
For any $u<s$,
\begin{multline}\label{eq:varrho}
\varrho_2(u,s)\mathrm{d}u\mathrm{d}s=
\E\left(\mathrm{d}X_u\mathrm{d}X_s\right)=
\E\left(\mathrm{d}X_u\E\left(\mathrm{d}X_s|\mathcal{F}_s\right)\right)\\
=\E\left(\mathrm{d}X_uw(X_s)\right)\mathrm{d}s \le
D\E\left(\mathrm{d}X_u\E\left(X_s|\mathcal{F}_{u^+}\right)\right)\mathrm{d}s\\
\le D\E\left(\mathrm{d}X_u {f}_r^{-1}\left(\E\left(f_r(X_s)|\mathcal{F}_{u^+}\right)\right)\right)\mathrm{d}s=
D\E\left(\mathrm{d}X_u {f}_r^{-1}\left(e^{r(s-u)}f_r(X_{u^+})\right)\right)\mathrm{d}s.
\end{multline}
Here we have used assumption (\ref{eq:felt1}), Jensen's inequality ($f_r$ is convex), and
also (\ref{eq:martingale}).

In order to estimate (\ref{eq:varrho}), we need appropriate upper and lower bounds for $f_r(k)$. 
There are two cases that we treat differently, depending on the value of $\alpha$ appearing
in (\ref{eq:felt2}):

\emph{(a)} $\alpha =1$ and 

\emph{(b)} $0< \alpha < 1$ .

Let us begin with case \emph{(a)}, when $w(k)=k+v(k)$ with
$v(k)=o(k)$ and $\sum _{k=1}^{\infty} \frac{v(k)}{k^2}< \infty$.

Throughout the following estimates the only thing we use is that for $x>0$, 
$x-\frac{x^2}{2}\le \log (1+x)\le x$. Recalling (\ref{eq:frdef}) that defines $f_r$, and
using $D$ again for denoting every constant,

\begin{multline*}
\log f_r(k)= 
\sum_{j=0}^{k-1}\log \left(1+\frac{r}{j+v(j)}\right)\le 
r \sum_{j=1}^{k-1}\frac{1}{j+v(j)}\\
=\frac{r}{v(0)}+ r \sum_{j=1}^{k-1}\left( \frac 1j -\frac{v(j)}{j(j+v(j))}\right)\le
r\log k + D,
\end{multline*}
using that $\frac{v(j)}{j(j+v(j))}$ is summable.

This way we obtain
\[f_r(k)\le D k^r.\]

Similarly
\[
\log f_r(k)\ge 
r \sum_{j=0}^{k-1}\frac{1}{j+v(j)}-r^2\sum_{j=0}^{k-1}\frac{1}{(j+v(j))^2}\ge
r\log k + D,
\]
with another finite $D$, so
\[
f_r(k)\ge D k^r,
\]
and this gives the upper bound for $f_r^{-1}$:
\[
f_r^{-1}(t)\le D t^{1/r}.
\]

This way we can continue (\ref{eq:varrho}) to obtain
\begin{multline}\label{eq:estim}
\varrho_2(u,s)\mathrm{d}u\mathrm{d}s\le 
D \E\left(\mathrm{d}X_u f_r^{-1}\left(e^{r(s-u)}X_{u^+}^r\right)\right)\mathrm{d}s\le
D \E\left( \mathrm{d}X_u e^{(s-u)}X_{u^+} \right)\mathrm{d}s\\
=D \E\left( w(X_u) e^{(s-u)}(X_u+1) \right)\mathrm{d}u \mathrm{d}s=
D e^{(s-u)} \E\left( X_u^2+X_u \right)\mathrm{d}u \mathrm{d}s.
\end{multline}

We now need estimates of $\E(X_u)$ and $\E(X_u^2)$.
The first one goes exactly like we have seen before: with $r$ large enough for ${f}_r$
to be convex,
\[
\E X_u \le f_r^{-1} \left( \E f_r (X_u) \right) = f_r^{-1}(e^{ru}) \le D e^u.
\]

To estimate $\E(X_u^2)$, consider the following. For arbitrary $r$ let us define
a function $g_r$ first at special coordinates: $g_r(k^2):=f_r(k)$ for $k\in \N$. Then for other 
$x\in \R^+$ let $g_r$ be extended piecewise linearly. This extension is convex if for all $k$
\[
\frac{g(k^2)-g((k-1)^2)}{k^2-(k-1)^2}\le \frac{g((k+1)^2)-g(k^2)}{(k+1)^2-k^2},
\]
which is equivalent with 
\begin{equation}\label{eq:ak}
r\ge \frac{2k+1}{2k-1}w(k)-w(k-1)=w(k)-w(k-1)+\frac{2}{2k-1}w(k).
\end{equation}
Since we have assumption (\ref{eq:felt1}) and $w(k)=k+o(k)$, the right-hand side of (\ref{eq:ak})
is bounded. This way we can find another (possibly larger) $r<\infty$ for which the extension
of $g_r$ is convex. We fix such an $r$ and continue:

\begin{equation}\label{eq:estimx2}
\E(X_u^2)\le g_r^{-1}(\E g_r(X_u^2))=g_r^{-1}(\E f_r(X_u))=g_r^{-1}(e^{ru}).
\end{equation}
Since
\[
g_r(x)\ge f_r(\lfloor \sqrt{x}\rfloor)\ge D x^{r/2}
\]
we obtain
\[
g_r^{-1}(t)\le Dt^{2/r}.
\] 
Therefore we can see from (\ref{eq:estimx2}) that
\[
\E(X_u^2)\le D e^{2u}.
\]

Let us now turn back to (\ref{eq:estim}) and observe that
\begin{multline*}
\int_u^{\infty}\int_0^{\infty}e^{-\lambda^*(u+s)}\varrho_2(u,s)\mathrm{d}u\mathrm{d}\\ 
\le 
D\int_u^{\infty}\int_0^{\infty}e^{-\lambda^*(u+s)}e^{s-u}e^{2u}\mathrm{d}u\mathrm{d}s + 
D\int_u^{\infty}\int_0^{\infty}e^{-\lambda^*(u+s)}e^{s-u}e^u\mathrm{d}u\mathrm{d}s.
\end{multline*}

In the exponents we have
\[
-\lambda^*u+(1-\lambda^*)s
\]
and
\[
(1-\lambda^*)u+(1-\lambda^*)s.
\]
Now recall that in the case $\alpha =1$, $\hat{\varrho}(1)=\infty$  
so $\lambda^*>1$, therefore the above integral is finite.

With this the proof of the Lemma is complete in the case $\alpha =1$.

Now let us turn to case \emph{(b)}, when $0<\alpha<1$, and $w(j)=j^{\alpha}+v(j)$ with
$v(j)=o(j)$. The estimation of $f_r(k)$: 

\begin{multline*}
\log f_r(k)= 
\sum_{j=0}^{k-1}\log \left(1+\frac{r}{j^{\alpha}+v(j)}\right)
\le 
\frac{r}{v(0)} + r \sum_{j=1}^{k-1}\left( \frac {1}{j^{\alpha}} -
\frac{v(j)}{j^{\alpha}(j^{\alpha}+v(j))}\right).
\end{multline*}
This time notice that
\[
\left| 
\sum_{j=0}^{k-1} \frac{v(j)}{j^{\alpha}(j^{\alpha}+v(j))} 
\right|\le
D+D\sum_{j=0}^{k-1}\frac{|v(j)|}{j^{\alpha}}\frac{1}{j^{\alpha}}\le
D+Dk^{1-\alpha}
\]
since $v(j)=o(j)$.
Therefore
\[
\log f_r(k)\le D+Dk^{1-\alpha}
\]
so
\[
f_r(k)\le De^{Dk^{1-\alpha}}.
\]
For the estimate from below consider that
\[
\log f_r(k)\ge
D\sum \left[  \frac{1}{j^{\alpha}+v(j)}-\frac{1}{(j^{\alpha}+v(j))^2}\right], 
\]
the absolute value of the sum of the second terms being of order $O(k^{1-2\alpha})$, and this
way we see that
\[
f_r(k)\ge De^{Dk^{1-\alpha}},
\]
so
\[
f_r^{-1}(e^{ru})\le Du^{\frac{1}{1-\alpha}}.
\]

Using these estimates we continue (\ref{eq:varrho}) with
\begin{multline*}
\varrho _2(u,s)\mathrm{d}u\mathrm{d}s\le
D\E \left( 
\mathrm{d}X_u (s-u+D+DX_{u^+}^{1-\alpha})^{\frac{1}{1-\alpha}} 
\right)\mathrm{d}s\\
\le D(s-u)^{\frac{1}{1-\alpha}}\E(X_{u^+}\mathrm{d}X_u)\mathrm{d}s\le
D(s-u)^{\frac{1}{1-\alpha}}\E(X_u^2+X_u)\mathrm{d}s.
\end{multline*}
The estimates of $\E(X_u)$ and $\E(X_u^2)$ go the same way as in case \emph{(a)}, we get
$\E(X_u)\le Du^{\frac{1}{1-\alpha}}$ and $\E(X_u^2)\le Du^{\frac{2}{1-\alpha}}$,
so the integral is obviously finite, since $\varrho_2(u,s)$ 
grows at most polynomially in $u$ and $s$.

This completes the proof of the Lemma.

\end{prf}

\section{Proof of the theorem}\label{chap:threesteps}

\subsection{First step: Characterization of ${n}_k(t)$}\label{chap:n_k}

Recall that $N_k(t)$ denotes the number of vertices that have at least
$k$ children at 
time $t$, its expectation is denoted by
$n_k(t):=\E\left(N_k(t)\right)$. We define $n_k(s):=0$ for $s<0$.
Similarly let $N_k^{(j)}(t)$ denote the same quantity but only in the subgraph which is
"under" the $j^{\mathrm{th}}$ child of the root (we define 
it to be 0 for $t<\tau_j$). The usefulness of this notation is that given 
$\tau:=\left(\tau _j\right)_{j=1}^{\infty}$, $N_k^{(j)}(t)$ and $N_k(t-\tau_j)$ have the same 
conditional distribution, so
\begin{equation}\label{eq:exp}
\E\left(N_k^{(j)}(t)\mid \tau\right)=\E\left(N_k(t-\tau _j)\mid \tau\right)=
n_k(t-\tau_j).
\end{equation}

Observe that the following simple recursion holds:
\begin{equation}
N_k(t)=\ind_{\{\tau_k<t\}}+\sum_{j=1}^{\infty} N_k^{(j)}(t). \label{eq:rekurzio}
\end{equation}
Note that for all $t\in (0,\infty)$, the sum above is almost surely finite.

Let us take the expectation of both sides of (\ref{eq:rekurzio})
in two steps: first taking conditional expectation
given $\tau$,
\[\E \left(N_k(t)\mid \tau\right)=\ind_{\{\tau_k<t\}}+
\sum_{j=1}^{\infty}n_k(t-\tau _j)=
\ind_{\{\tau_k<t\}}+\int_0^t n_k(t-s)\mathrm{d}X_s.\]
Here we used (\ref{eq:exp}).

Then taking the expectation of both sides we get

\begin{equation}\label{eq:nkt}
n_k(t)=\pr(\tau _k<t)+\int_0^t n_k(t-s)\varrho (s)\mathrm{d}s.
\end{equation}

Taking the Laplace transform of both sides
and defining $p_k(t):=\pr(\tau _k<t)$ we get
\[\hat{n}_k(\lambda)=\hat{p}_k(\lambda)+\hat{n}_k(\lambda)\hat{\varrho}(\lambda).\]

We have seen in the previous chapter that $(1-\hat{\varrho}(\lambda))$ is 
positive for every $\lambda>\lambda ^*$, so for these $u$
\[\hat{n}_k(\lambda)=\frac{\hat{p}_k(\lambda)}{1-\hat{\varrho}(\lambda)}.\]

Regarding $\frac{\hat{p}_k(\lambda)}{1-\hat{\varrho}(\lambda)}$ as a complex function, it has a simple
pole at $\lambda=\lambda^*$ (since $\lambda\mapsto \hat{\varrho}(\lambda)$ strictly dicreases on the 
interval where it is 
finite, $\lambda ^*$ is a simple root of $1-\hat{\varrho}(\lambda)$).
This way $\hat{n}_k(\lambda)$ has the form
\[
\hat{n}_k(\lambda)=\frac{\hat{p}_k(\lambda)}{1-\hat{\varrho}(\lambda)}=
\frac{a\hat{p}_k(\lambda^*)}{\lambda-\lambda^*}+\hat{m}_k(\lambda),
\]
where $\hat{m}_k(\lambda)$ is analytic at $\lambda=\lambda ^*$. This way $\hat{m}_k(\lambda)$ is the Laplace 
transform of some function $m_k(t)$ for which $\lim_{t\to \infty}e^{-\lambda ^* t}m_k(t)=0$. 
Therefore we have
\[n_k(t)=a\hat{p}_k(\lambda^*)e^{\lambda^{^*} t}+m_k(t)\] 
and we can conclude that
\[\lim_{t\to \infty}e^{-\lambda ^* t}n_k(t)=a\hat{p}_k(\lambda ^*)=:d_k.\]

We define $c_k$ to be the fraction of $d_k$ and $d_0$,
\[\frac{d_k}{d_0}=:c_k.\]

Since $n_k(t)$ is of order $e^{\lambda ^*t}$, it is useful to
introduce the notation $\tilde{n}_k(t):=e^{-\lambda ^* t}n_k(t)$.

In order to compute $c_k$ explicitly, first recall (\ref{eq:varrhokalap}), then consider that
we can also get $\hat{p}_k(\lambda ^*)$ simply: for $k\ge 1$
\[\hat{p}_k(\lambda)=\int_0^{\infty}e^{-\lambda t}\pr (\tau _k<t)\mathrm{d}t=
\frac{1}{\lambda}\E (e^{-\lambda \tau_k})=\frac{1}{\lambda}\prod_{j=0}^{k-1}\frac{w(j)}{\lambda+w(j)}.\]

This gives us the exact values of the $c_k$:
\[c_k=\frac{d_k}{d_0}=\prod_{j=0}^{k-1}\frac{w(j)}{\lambda ^*+w(j)}\]
since $\hat{p}_0(\lambda ^*)=\frac{1}{\lambda ^*}$.

\subsection{Second step: Characterization of $n_{k,l}(t)$}
\label{chap:n_{k,l}}

Let us define $n_{k,l}(t):=\E (N_k(t)N_l(t))$, and introduce the notation 
$\tilde{n}_{k,l}(t):=e^{-2\lambda ^* t}n_{k,l}(t)$. The aim of this section is to prove
that
\[
\lim_{t \to \infty} \tilde{n}_{k,l}(t)=Cd_kd_l
\] 
for some constant $0<C<\infty$.

Using the recursion given by (\ref{eq:rekurzio}) we can see that, for $k\ge l$,
\begin{multline*}
N_k(t)N_l(t)=
\ind_{\{\tau_k<t\}}+
\ind_{\{\tau_k<t\}}\sum_{j=1}^{\infty}N^{(j)}_l(t)+
\ind_{\{\tau_l<t\}}\sum_{j=1}^{\infty}N^{(j)}_k(t)\\
+\sum_{j=1}^{\infty}N^{(j)}_k(t)N^{(j)}_l(t)+
\sum_{i\ne j}^{\infty}N^{(i)}_k(t)N^{(j)}_l(t).
\end{multline*}
Taking conditional expectation given $\tau$, then expectation with respect to
$\tau$ (as done before), we get
\begin{multline}\label{eq:thelimit}
n_{k,l}(t)=
\pr (\tau_k<t )\\
+\E \left(\ind_{\{\tau_l<t\}}\int_0^tn_k(t-s)\mathrm{d}X_s+
\ind_{\{\tau_k<t\}}\int_0^tn_l(t-s)\mathrm{d}X_s\right)\\
+\int_0^tn_{k,l}(t-s)\varrho (s)\mathrm{d}s+
\int_0^t\int_0^tn_k(t-s)n_l(t-u)\varrho_2(u,s)\mathrm{d}u\mathrm{d}s.
\end{multline}

After multiplying the equation by $e^{-2\lambda ^*t}$, we can easily identify the limit 
(as $t \to \infty$) of the first, second and fourth terms of the right hand side:

Since $\pr(\tau _k<t)<1$, the limit of the first term divided by $e^{2\lambda^* t}$ is trivially 0. 
For the second term consider that
\begin{multline*}
\E \left(\ind_{\{\tau_l<t\}}\int_0^te^{-2\lambda ^*t}n_k(t-s)\mathrm{d}X_s\right)\le 
e^{-2\lambda ^* t}\E \left(\int_0^tn_k(t-s)\mathrm{d}X_s\right)\\
=e^{-2\lambda ^* t}\int_0^tn_k(t-s)\varrho (s)\mathrm{d}s=e^{-2\lambda ^* t}(n_k(t)-p_k(t))
\end{multline*} 
from (\ref{eq:nkt}),
and we know that $n_k(t)$ is of order $e^{\lambda ^* t}$, and $p_k(t)$ is bounded. So the
second term in (\ref{eq:thelimit}) divided by $e^{2\lambda ^* t}$ is also 0 in the limit.

As for the fourth term of (\ref{eq:thelimit}), using the statement and notations of 
the Lemma (see (\ref{eq:kell})), we can conclude that
\[
\lim_{t \to \infty}
\int_0^t\int_0^t\tilde{n}_k(t-s)\tilde{n}_l(t-u)e^{-\lambda ^*(u+s)}\varrho_2(u,s)
\mathrm{d}u\mathrm{d}s
=B_2d_kd_l
\]
by dominated convergence (we already know the convergence and therefore the boundedness of 
$\tilde{n}_k(t)$ and $\tilde{n}_l(t)$).

This way we see
\begin{equation}\label{eq:nklt}
\tilde{n}_{k,l}(t)=
\int_0^t\tilde{n}_{k,l}(t-s)e^{-2\lambda ^* s}\varrho(s)\mathrm{d}s+
B_2d_kd_l+
\varepsilon_{k,l}(t),
\end{equation}
where $\varepsilon_{k,l}(t)\to 0$ as $t\to \infty$.

Now let us assume for a moment that the limit $d_{k,l}:=\lim_{t \to \infty} \tilde{n}_{k,l}(t)$
does exist. In this case dominated convergence could also be used in the third term of
(\ref{eq:thelimit}), and the following would be true:
\[
\lim_{t \to \infty}\int_0^te^{-2\lambda ^* t}n_{k,l}(t-s)\varrho (s)\mathrm{d}s=
\lim_{t \to \infty}\int_0^t\tilde{n}_{k,l}(t-s)e^{-2\lambda ^* s}\varrho (s)\mathrm{d}s=
B_1d_{k,l},
\]
recall the notation in (\ref{eq:b1def}), and that $B_1<1$. 

This way \textbf{if} $\tilde{n}_{k,l}(t)$ was convergent, then its limit could only be
\[
d_{k,l}=\frac{B_2}{1-B_1}d_kd_l.
\]

To show that the limit really exists, first note that $\tilde{n}_{k,l}(t)$ is 
bounded: for $M_{k,l}(t):=\sup_{s<t}\tilde{n}_{k,l}(s)$ we 
get that
\[M_{k,l}(t)\le E_{k,l}+M_{k,l}(t)B_1+M_kM_lB_2,\]
and so $M_{k,l}(t)$ is bounded by a constant independent of $t$. (Here $E_{k,l}$ is an upper 
bound for $\varepsilon_{k,l}(t)$,
$M_k=\sup_{s>0}\tilde{n}_k(s)$ and the $B_i$ are as before.) 

Now as for the convergence, let us introduce 
$\lev{n}_{k,l}(t):=\tilde{n}_{k,l}(t)-\frac{B_2}{1-B_1}d_kd_l$ and rearrange 
equation (\ref{eq:nklt}),

\[
\lev{n}_{k,l}(t)=
\lev{\varepsilon}_{k,l}(t)+
\int_0^t\lev{n}_{k,l}(t-s)e^{-2\lambda^*s}\varrho(s)\mathrm{d}s,
\]
where $\lev{\varepsilon}_{k,l}(t)\to 0$ as $t\to \infty$. 

Since we have shown that $\tilde{n}_{k,l}(t)$ is bounded, so is $\lev{n}_{k,l}(t)$. Let
$\lev{M}_{k,l}(t):=\sup_{s\ge t}|\lev{n}_{k,l}(s)|$, 
$\lev{E}_{k,l}(t):=\sup_{s\ge t}|\lev{\varepsilon}_{k,l}(s)|$,
and fix arbitrarily $0<u<t_0$. 
For these and for all $t>t_0$ 
\begin{multline*}
|\lev{n}_{k,l}(t)|\le |\lev{\varepsilon}_{k,l}(t)|+
\left|\int_0^u\lev{n}_{k,l}(s)e^{-2\lambda ^*(t-s)}\varrho (t-s)\mathrm{d}s\right|\\
+\left|\int_u^t\lev{n}_{k,l}(s)e^{-2\lambda ^*(t-s)}\varrho (t-s)\mathrm{d}s\right|
\end{multline*}
so
\[
|\lev{n}_{k,l}(t)|\le \lev{E}_{k,l}(t_0)+
e^{-\lambda ^* (t-u)}\lev{M}_{k,l}(0)+
\lev{M}_{k,l}(u)B_1.
\]
This way
\[\lev{M}_{k,l}(t_0)\le \lev{E}_{k,l}(t_0)+e^{-\lambda ^* (t_0-u)}\lev{M}_{k,l}(0)+
\lev{M}_{k,l}(u)B_1.\]
Letting $t_0\to \infty$ with $u$ remaining fixed
\[\lev{M}_{k,l}(\infty)\le \lev{M}_{k,l}(u)B_1,\]
and now letting $u\to \infty$
\[\lev{M}_{k,l}(\infty)\le \lev{M}_{k,l}(\infty)B_1.\]
Since $B_1<1$ this means that $\lev{M}_{k,l}(\infty)=0$, so $\tilde{n}_{k,l}(t)$ is
convergent and $\lim_{t\to \infty}\tilde{n}_{k,l}(t)=\frac{B_2}{1-B_1}d_kd_l$.

\subsection{Third step: End of proof}

We have shown that for all $k,l$
\[\E\left((d_l\tilde{N}_k(t)-d_k\tilde{N}_l(t))^2\right)\to 0\quad \text{as }t\to \infty,\]
the random variable $d_l\tilde{N_k}(t)-d_k\tilde{N}_l(t)$ converges to 0 in
$L_2$, so in probability as well.
From this we can argue that for any fixed $\eta>0$ and $\delta>0$
\begin{multline*}
\pr \left(\left|\frac{N_k(t)}{N_0(t)}-\frac{d_k}{d_0}\right|>\delta\right)=
\pr \left(\left|\frac{N_k(t)}{N_0(t)}-\frac{d_k}{d_0}\right|>\delta \biggm|
\tilde{N}_0(t)\ge \eta \right)\pr (\tilde{N}_0(t)\ge \eta)\\
+\pr \left(\left|\frac{N_k(t)}{N_0(t)}-\frac{d_k}{d_0}\right|>\delta \biggm|
\tilde{N}_0(t)<\eta \right)\pr (\tilde{N}_0(t)<\eta)\\
\le 
\pr (|d_0\tilde{N}_k(t)-d_k\tilde{N}_0(t)|>d_0\eta \delta)+\pr (\tilde{N}_0(t)<\eta).
\end{multline*}

Now with fixed $\eta$ let $t\to \infty$ to see that 
\begin{equation}\label{eq:ineq}
\varlimsup_{t\to \infty}\pr \left(\left|\frac{N_k(t)}{N_0(t)}-\frac{d_k}{d_0}\right|>
\delta \right)\le
\varlimsup_{t\to \infty}\pr (\tilde{N}_0(t)<\eta).
\end{equation}

Therefore it is sufficient to show that in the limit $t \to \infty$ the random variable 
$\tilde{N}_0(t)$ does not have a positive mass at 0.

Consider that
\[
N_0(t)=1+\sum_{j=1}^{\infty} N_0^{(j)}(t)\stackrel{\textrm{d}}{=}
1+\sum_{j=1}^{\infty} \ind _{(\tau _j<t)}N_{0,j}(t-\tau _j),
\]
where $\left(N_{0,j}\right)_{j>0}$ are independent identically distributed random processes, 
with the same distribution as $N_0$. So
\begin{equation}\label{eq:n0rek}
\frac{N_0(t)}{e^{\lambda^*t}}\stackrel{\textrm{d}}{=}
e^{-\lambda^*t}+\sum_{j=1}^{\infty} \ind _{(\tau _j<t)}e^{-\lambda^*\tau _j}
\frac{N_{0,j}(t-\tau _j)}{e^{\lambda^*(t-\tau _j)}}.
\end{equation}

Since $\left(\frac{N_0(t)}{e^{\lambda^*t}}\right)_{t\ge 0}$ is tight, in every subsequence 
there is a sub-subsequence $(t_n)_{n\ge 0}$ along which $\frac{N_0(t)}{e^{\lambda^*t}}$ 
converges weakly to some random variable $\xi$. By (\ref{eq:n0rek}) for this variable
\[
\xi\stackrel{\textrm{d}}{=}\sum_{j=1}^{\infty}e^{-\lambda^*\tau _j}\xi _j,
\]
where the $\xi _j$ are iid with the same distribution as $\xi $.

This means that
\[
\pr (\xi =0)=\pr (\xi _j=0 \text{ for all }j)=
\lim _{k\to {\infty}}\left(\pr (\xi =0)\right)^k.
\]
It follows that if $\xi$ had a positive mass at 0, then $\xi$ would be a random variable that is 
almost surely 0. Since 
$\E\left(\frac{N_0(t)}{e^{\lambda^*t}}\right) \to d_0>0$, this could only happen if
$\E\left(\frac{N_0^2(t)}{e^{2\lambda^*t}}\right) \to \infty$, but we have already shown that
this is not true, it converges to a finite limit $Cd_0^2$.

We can finally conclude that the right hand side of (\ref{eq:ineq}) converges to 0 as
$\eta \to 0$, so the theorem holds:
\[\frac{N_k(t)}{N_0(t)}\to \frac{d_k}{d_0}=c_k\quad\text{in probability as }t\to \infty .\] 


The author thanks B\'alint T\'oth for valuable discussions.

\vfill\vfill\vfill

\hbox{\sc
\vbox{\noindent
\hsize75mm
Anna Rudas\\
Institute of Mathematics\\
Technical University Budapest\\
Egry J\'ozsef u. 1.\\
H-1111 Budapest, Hungary\\
{\tt rudasa{@}math.bme.hu}
}
}


\begin{thebibliography}{9}
\bibitem{Barabasi}Barab\'asi, A.-L., Albert, R., Emergence of Scaling in Random
Networks, \textit{Science} \textbf{286} (1999), 509-512
\bibitem{Bollobas_summary}Bollob\'as, B., Mathematical results on scale-free random graphs, in: 
S.~Bornholdt and H.~G.~Schuster, editors, \textit{Handbook of Graphs and Networks}, pages 1-34.
Wiley, 2002.
\bibitem{Bollobas_degree}Bollob\'as, B., Riordan, O., Spencer, J., Tusnády, G., The Degree Sequence of a 
Scale-free Random Graph Process, \textit{Random Structures and Algorithms} \textbf{18} (2001), 279-290. 
\bibitem{Mori}Móri, T., On Random Trees, \textit{Studia Sci. Math. Hungar.}, \textbf{39} (2002), 
143-155.
\bibitem{Szymanski}Szyma\'nski, J., On a Nonuniform Random Recursive Tree, 
\textit{Ann. Discrete Math.} \textbf{33} (1987), 297-306.
\end{thebibliography}
\end{document}